\documentclass{amsart}
\usepackage{graphicx}
\usepackage{psfrag}

\theoremstyle{plain}
\newtheorem{prop}{Proposition}
\newtheorem{thm}{Theorem}
\newtheorem{cor}{Corollary}
\newtheorem{lemma}{Lemma}
\newtheorem{conj}{Conjecture}

\theoremstyle{remark}
\theoremstyle{definition}
\newtheorem{remark}{Remark}

\newcommand{\ot}{\leftarrow}
\newcommand{\Schub}{{\mathfrak S}}
\newcommand{\Groth}{{\mathfrak G}}

\newcommand{\Gr}{\operatorname{Gr}}

\newcommand{\Z}{{\mathbb Z}}

\newcommand{\C}{{\mathbb C}}

\newcommand{\bull}{{\scriptscriptstyle \bullet}}
\DeclareMathOperator{\Hom}{Hom}

 \DeclareMathOperator{\codim}{codim}

\DeclareMathOperator{\GL}{GL} \DeclareMathOperator{\Ext}{Ext}

\newcommand{\ov}{\overline}

\newcommand{\wt}{\widetilde}

\newcommand{\bu}{{\mathbf u}}
\newcommand{\bw}{{\mathbf w}}

\newcommand{\euler}{{\mathcal E}}

\newcommand{\pic}[2]{\includegraphics[scale=0.#1]{#2}}

\newcommand{\ignore}[1]{}

\psfrag{i}{$i$}
\psfrag{ip1}{$i+1$}
\psfrag{jm1}{$j-1$}
\psfrag{j}{$j$}

\begin{document}

\title[A formula for non-equioriented quiver orbits of type $A$]
{A formula for non-equioriented \\ quiver orbits of type $A$}
\author{Anders Skovsted Buch}
\address{Department of Mathematics, Rutgers University, 110
  Frelinghuysen Road, Piscataway, NJ 08854, USA}
\email{asbuch@math.rutgers.edu}
\author{Rich\'ard Rim\'anyi}
\address{Department of Mathematics, The University of North Carolina
  at Chapel Hill,\linebreak CB \#3250, Phillips Hall, Chapel Hill, NC
  27599, USA}
\email{rimanyi@email.unc.edu}

\date{\today}

\begin{abstract}
  We prove a positive combinatorial formula for the equivariant class
  of an orbit closure in the space of representations of an arbitrary
  quiver of type $A$.  Our formula expresses this class as a sum of
  products of Schubert polynomials indexed by a generalization of the
  minimal lace diagrams of Knutson, Miller, and Shimozono.  The proof
  is based on the interpolation method of Feh\'er and Rim\'anyi.  We
  also conjecture a more general formula for the equivariant
  Grothendieck class of an orbit closure.
\end{abstract}

\maketitle

\section{Introduction}

A {\em quiver\/} is an oriented graph $Q = (Q_0, Q_1)$ consisting of a
set of vertices $Q_0$ and a set of arrows $Q_1$.  Each arrow $a \in
Q_1$ has a tail $t(a) \in Q_0$ and a head $h(a) \in Q_0$.  In this
paper we will consider a quiver $Q$ of type $A$, i.e.\ a chain of
vertices with arrows between them.  We identify the vertex and arrow
sets with integer intervals, $Q_0 = \{0,1,2,\dots,n\}$ and $Q_1 =
\{1,2,\dots,n\}$, such that $\{t(a),h(a)\} = \{a-1,a\}$ for each $a
\in Q_1$.  We also set $\delta(a) = h(a)-t(a)$, which equals $-1$ for
a {\em leftward\/} arrow and $+1$ for a {\em rightward\/} arrow.

Fix a dimension vector $e = (e_0,e_1,\dots,e_n)$ of non-negative
integers, and set $E_i = \C^{e_i}$ for each $i$.  The set of quiver
representations of dimension vector $e$ form the affine space
\[ V = \Hom(E_{t(1)},E_{h(1)}) \oplus \dots \oplus
   \Hom(E_{t(n)},E_{h(n)}) \,,
\]
which has a natural action (with finitely many orbits) of the group
$G = \GL(E_0) \times \dots 
\times \GL(E_n)$ given by $(g_0,\dots,g_n).(\phi_1,\dots,\phi_n) =
(g_{h(1)} \phi_1 g_{t(1)}^{-1}, \dots, g_{h(n)} \phi_n
g_{t(n)}^{-1})$.  The goal of this paper is to prove a formula for the
$G$-equivariant cohomology class of an orbit closure for this action.
We note that Poincar\'e duality in equivariant cohomology was
introduced by Kazarian \cite{kazarian:characteristic}, but simpler
methods can be used to define the classes of Zariski closed subsets of
$V$ \cite{feher.rimanyi:calculation, fulton:notes}.  Our formula can
also be interpreted as a formula for degeneracy loci defined by a
quiver of vector bundles and bundle maps over a complex variety.  This
application relies on Bobi{\'n}ski and Zwara's proof that orbit
closures of type $A$ are Cohen-Macaulay
\cite{bobinski.zwara:normality}.

The quiver $Q$ is {\em equioriented\/} if all arrows have the same
direction.  A formula for the orbit closures for such a quiver was
proved by Buch and Fulton \cite{buch.fulton:chern}.  Notice also that
the problem specializes to the classical Thom-Porteous formula when
$n=1$.  The formula proved in this paper generalizes a different
formula for equioriented orbit closures, called the {\em component
  formula}, which was conjectured by Knutson, Miller, and Shimozono
and proved in \cite{knutson.miller.ea:four} and
\cite{buch.feher.ea:positivity}.

For an arbitrary quiver of Dynkin type, the {\em interpolation
  method\/} of Feh\'er and Rim\'anyi makes it possible to compute the
class of an orbit closure as the unique solution to a system of linear
equations, which say that this class must vanish when restricted to a
disjoint orbit \cite[\S2]{feher.rimanyi:calculation}.  The proof of
our formula relies on this method, as well as on a simplification of
the ideas from \cite{buch.feher.ea:positivity}.

The $G$-orbits in $V$ are classified by the {\em lace diagrams\/} of
Abeasis and Del Fra \cite{abeasis.del-fra:degenerations,
  abeasis.del-fra:degenerations*2}.  For equioriented quivers, these
diagrams were reinterpreted as sequences of permutations by Knutson,
Miller, and Shimozono \cite{knutson.miller.ea:four}, who called a lace
diagram {\em minimal\/} if the sum of the lengths of these
permutations is equal to the codimension of the corresponding orbit.
The component formula writes the class of an orbit closure as a sum of
products of Schubert polynomials indexed by all minimal lace diagrams
for the orbit.  The same construction turns out to work for an
arbitrary quiver of type $A$, although most definitions need to be
changed to take the orientation of the arrows into account, including
the definition of a minimal lace diagram.  By combining our definition
of non-equioriented minimal lace diagrams with certain $K$-theoretic
transformations of lace diagrams from \cite{buch.feher.ea:positivity},
we furthermore obtain a natural conjecture for the equivariant
Grothendieck class of an orbit closure.  This conjecture generalizes
the $K$-theoretic component formulas from \cite{buch:alternating,
  miller:alternating}.

Our paper is organized as follows.  In Section~\ref{S:formula} we give
the definition of minimal lace diagrams, state our formula, and prove
some combinatorial properties of the formula.  We also explain its
interpretation as a formula for degeneracy loci.
Section~\ref{S:proof} explains the interpolation method and completes
the proof of our formula.  Finally, in Section~\ref{S:ktheory} we pose
our conjectured formula for the Grothendieck class of an orbit closure
of type $A$.

We thank the referee for several helpful suggestions to our exposition.

\section{The non-equioriented component formula}
\label{S:formula}

\subsection{Lace diagrams} \label{S:lacediag}

The $G$-orbits in $V$ are classified by the lace diagrams of Abeasis
and Del Fra \cite{abeasis.del-fra:degenerations}.  Define a {\em lace
  diagram\/} for the dimension vector $e$ to be a sequence of $n+1$
columns of dots, with $e_i$ dots in column $i$, together with line
segments connecting dots of consecutive columns.  Each dot may be
connected to at most one dot in the column to the left of it, and to
at most one dot in the column to the right of it.

The quiver representations $\phi = (\phi_1,\dots,\phi_n)$ in the orbit
given by a lace diagram can be obtained by identifying the dots of
column $i$ with chosen basis vectors of $E_i$, and defining each
linear map $\phi_a : E_{t(a)} \to E_{h(a)}$ according to the
connections between the dots.  In other words, if dot $j$ of column
$t(a)$ is connected to dot $k$ of column $h(a)$, then $\phi_a$ maps
the $j$th basis element of $E_{t(a)}$ to the $k$th basis element of
$E_{h(a)}$; and if dot $j$ of column $t(a)$ is not connected to any
dot in column $h(a)$, then the corresponding basis element of
$E_{t(a)}$ is mapped to zero.  For example, the following lace diagram
represents an orbit in the space of representations of the quiver $Q =
(\circ \to \circ \ot \circ \to \circ \to \circ)$ of dimension vector
$e = (3, 4, 4, 3, 3)$.
\[ \pic{50}{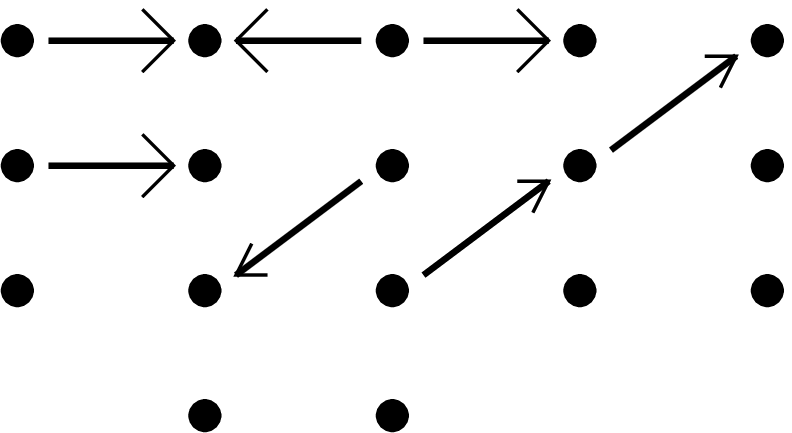} \]

A lace diagram can be interpreted as a sequence of $n$ permutations as
follows.  For each {\em rightward\/} arrow $a \in Q_1$ we let $w_a$ be
the permutation of smallest possible length such that $w_a(k) = j$
whenever the $k$th dot from the {\em top\/} of column $a$ is connected
to the $j$th dot from the {\em top\/} of column $a-1$.  If $a \in Q_1$
is a {\em leftward\/} arrow then we let $w_a$ be the permutation of
smallest length such that $w_a(j) = k$ if the $j$th dot from the {\em
  bottom\/} of column $a-1$ is connected to the $k$th dot from the
{\em bottom\/} of column $a$.  Notice in particular that each
permutation $w_a$ is read off the diagram {\em against\/} the
direction of the arrow $a \in Q_1$.  The lace diagram is determined by
the sequence of permutations $\bw = (w_1,\dots,w_n)$ together with the
dimension vector $e$.  Equivalently, the permutation sequence $\bw$
describes the connections between the dots of an {\em extension\/} of
the lace diagram, which is obtained by adding extra dots and
connections below each rightward arrow and above each leftward arrow.
The above displayed lace diagram corresponds to the permutation
sequence $(w_1,w_2,w_3,w_4)$ where $w_1 = 12453$, $w_2 = 536412$, $w_3
= 13524$, and $w_4 = 24513$.  The diagram has the following extension.
\[ \pic{50}{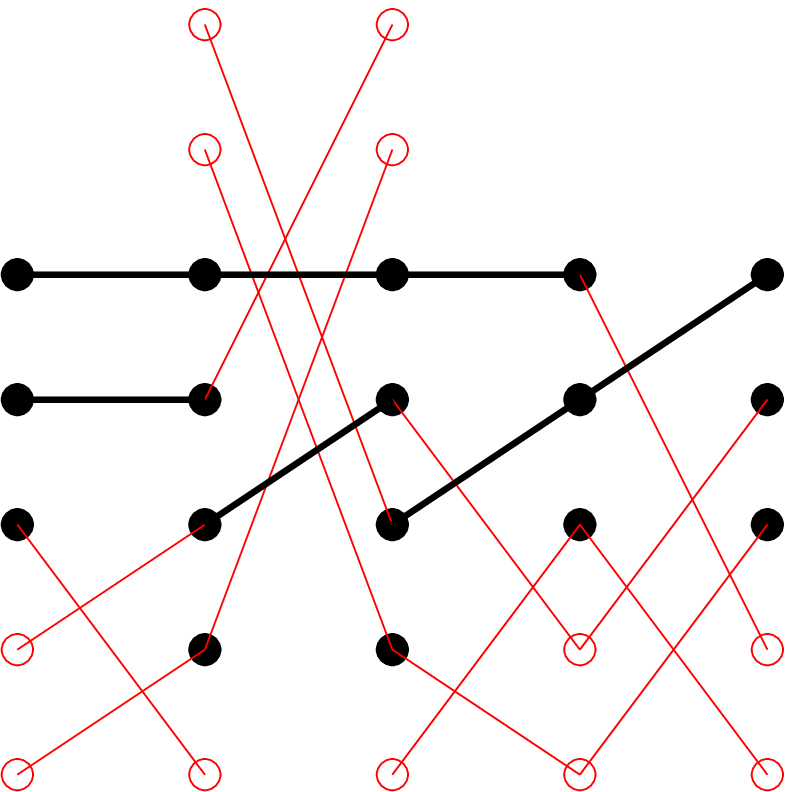} \]

A permutation $w$ is called a {\em partial permutation\/} from $p$
elements to $q$ elements if all descent positions of $w$ are smaller
than or equal to $p$, while the descent positions of $w^{-1}$ are
smaller than or equal to $q$.  In other words we have $w(i)<w(i+1)$
for $p > i$ and $w^{-1}(i) < w^{-1}(i+1)$ for $i > q$.  A sequence
$\bw = (w_1,\dots,w_n)$ of permutations represents a lace diagram if
and only if each permutation $w_a$ is a partial permutation from
$e_{h(a)}$ elements to $e_{t(a)}$ elements.  In the following we
identify a lace diagram with its permutation sequence $\bw$.

\subsection{Minimal lace diagrams}

A {\em strand\/} of a lace diagram is a maximal sequence of connected
dots and line segments, and the {\em extension\/} of a strand is
obtained by also including the extra line segments that it is directly
connected to in the extended lace diagram.  The {\em length\/} of the
lace diagram $\bw = (w_1,\dots,w_n)$ is the sum $\sum \ell(w_a)$ of
the lengths of the permutations $w_a$.  Equivalently, the length is
equal to the total number of crossings in the extended diagram of
$\bw$.

For an orbit $\mu \subset V$ and vertices $0\leq i \leq j \leq n$, we
define $s_{ij} = s_{ij}(\mu)$ to be the number of (non-extended)
strands starting at column $i$ and terminating at column $j$ for any
lace diagram representing $\mu$.  We also let $r_{ij} = r_{ij}(\mu)$
denote the total number of connections from column $i$ to column $j$,
i.e.\ $r_{ij} = \sum_{k\leq i, l \geq j} s_{kl}$.

\begin{lemma} \label{L:mincross}
The length of a lace diagram representing the orbit
  $\mu$ is greater than or equal to the number
\[ d(\mu) = \sum_{i<j:\,\delta(i+1)=\delta(j)}
   (r_{i+1,j}-r_{ij})(r_{i,j-1}-r_{ij}) \ +
   \sum_{i<j:\,\delta(i+1)\neq \delta(j)} r_{ij} s_{i+1,j-1} \,.
\]
\end{lemma}
\begin{proof}
  Consider vertices $i,j \in Q_0$ with $i<j$, and assume that the
  arrow between $i$ and $i+1$ has the same direction as the arrow
  between $j-1$ and $j$, that is $\delta(i+1) = \delta(j)$.  Since the
  left end of a strand starting at column $i+1$ is extended in the
  same direction (up or down) as the right end of a strand terminating
  at column $j-1$, it follows that (the extensions of) these strands
  must cross if the first strand passes through column $j$ and the
  second strand passes through column $i$.  There are exactly
  $(r_{i+1,j} - r_{ij})(r_{i,j-1} - r_{ij})$ examples of this.
\[ \pic{50}{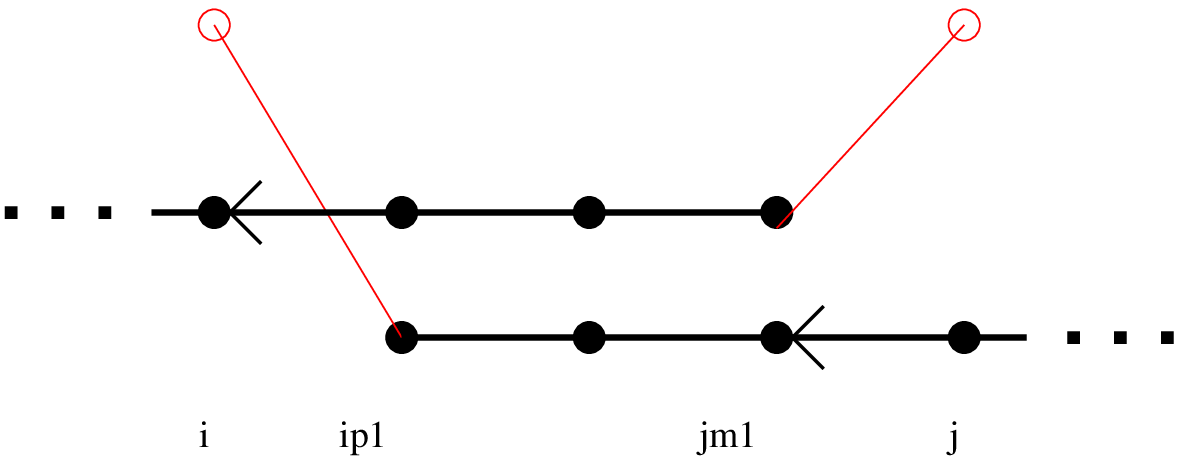} \]

On the other hand, if $\delta(i+1) \neq \delta(j)$, then the left and
right ends of a strand from column $i+1$ to column $j-1$ are extended
in opposite directions, which means that such a strand must cross all
strands connecting column $i$ to column $j$.  This happens in $r_{ij}
s_{i+1,j-1}$ examples.  We have therefore identified $d(\mu)$ forced
crossings in any lace diagram representing the orbit $\mu$.
\[ \pic{50}{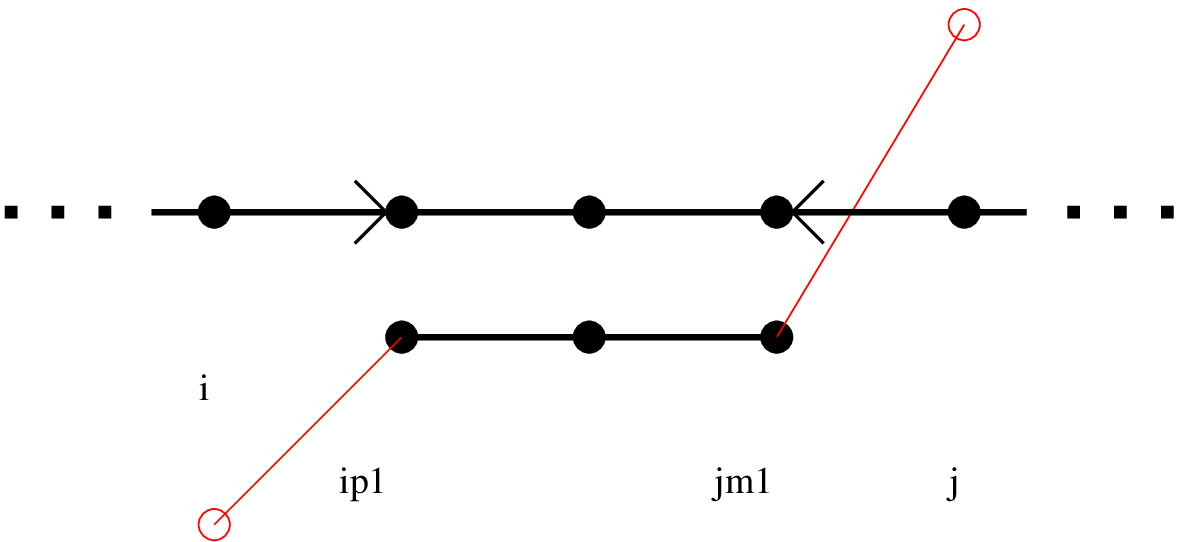} \]
\end{proof}

We will prove later that the integer $d(\mu)$ of
Lemma~\ref{L:mincross} is equal to the codimension of $\mu$ in $V$.
We will call a lace diagram for $\mu$ {\em minimal\/} if its length is
equal to $d(\mu)$.  This extends Knutson, Miller, and Shimozono's
definition of a minimal lace diagram for an equioriented quiver
\cite{knutson.miller.ea:four}.  The following extended lace diagram is
minimal and represents the same orbit as the diagrams of
Section~\ref{S:lacediag}.
\[ \pic{50}{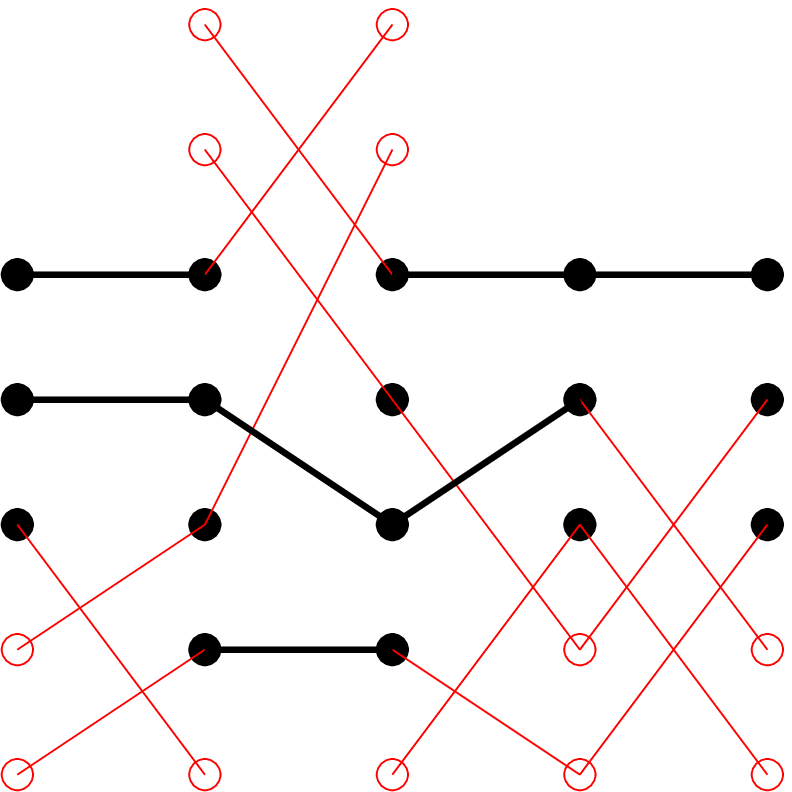} \]

Notice that a lace diagram is minimal if and only if any two strands
cross at most once, and not at all if they start or terminate at the
same column (cf.~\cite[Thm.~3.8]{knutson.miller.ea:four}).  In fact,
none of the forced crossings identified in the proof of
Lemma~\ref{L:mincross} involve strands starting or terminating in the
same column, and if two strands starting and terminating in different
columns are not forced to cross, then they cross an even number of
times.

\subsection{Schubert polynomials}

To state our formula, we need the Schubert polynomials of Lascoux and
Sch\"utzen\-berger \cite{lascoux.schutzenberger:polynomes}.  The {\em
  divided difference operator\/} $\partial_{a,b}$ with respect to two
variables $a$ and $b$ is defined by
\[ \partial_{a,b}(f) = \frac{f(a,b) - f(b,a)}{a-b} \]
where $f$ is any polynomial in these (and possibly other) variables.
The {\em double Schubert
  polynomials\/} $\Schub_w(x;y) =
\Schub_w(x_1,\dots,x_m;y_1,\dots,y_m)$ given by permutations $w \in
S_m$ are uniquely determined by the identity
\begin{equation} \label{E:schubdef}
  \partial_{x_i,x_{i+1}}(\Schub_w(x;y)) =
  \begin{cases} \Schub_{w s_i}(x;y) & \text{if $\ell(w s_i) < \ell(w)$}\\
    0 & \text{if $\ell(w s_i) > \ell(w)$} \end{cases}
\end{equation}
together with the expression $\Schub_{w_0}(x;y) = \prod_{i+j\leq m}
(x_i-y_j)$ for the longest permutation $w_0$ in $S_m$.  Using that
$\Schub_w(x;y) = (-1)^{\ell(w)} \Schub_{w^{-1}}(y;x)$, the identity
(\ref{E:schubdef}) is equivalent to
\begin{equation} \label{E:dualdef}
  \partial_{y_i,y_{i+1}}(\Schub_w(x;y)) =
  \begin{cases} - \Schub_{s_i w}(x;y) & \text{if $\ell(s_i w) < \ell(w)$}\\
    0 & \text{if $\ell(s_i w) > \ell(w)$}\,. \end{cases}
\end{equation}

For any permutations $u,w \in S_m$, the definition of Schubert
polynomials implies that the specialization $\Schub_w(y_u;y) =
\Schub_w(y_{u(1)},\dots,y_{u(m)}; y_1,\dots,y_m)$ is zero unless $w
\leq u$ in the Bruhat order on $S_m$, and for $u = w$ we have
\begin{equation} \label{E:spec}
  \Schub_u(y_{u(1)},\dots,y_{u(m)}; y_1,\dots,y_m) =
  \prod_{i<j:\,u(i)>u(j)} (y_{u(i)} - y_{u(j)}) \,.
\end{equation}
Furthermore, if $k$ and $l$ are the last descent positions of $w$ and
$w^{-1}$, respectively, then only the variables
$x_1,\dots,x_k,y_1,\dots,y_l$ occur in $\Schub_w(x;y)$.

\subsection{Statement of the formula}

For each $i \in Q_0$ we let $x^i = \{x^i_1,\dots,x^i_{e_i}\}$ be a set
of $e_i$ variables.  These variables are identified with the Chern
roots in $H^*_T(V)$ of the $i$th factor of $G$, where $T$ is a maximal
torus of $G$.  Then $H^*_T(V)$ is the polynomial ring $\Z[x^i_j \mid
0\leq i\leq n, 1\leq j\leq e_i]$ in these variables, and $H^*_G(V)
\subset H^*_T(V)$ is the subring of polynomials which are separately
symmetric in each set of variables $x^i$.  We let $\wt x^i = \{
x^i_{e_i}, \dots, x^i_1 \}$ denote the variables $x^i$ in the opposite
order.  Given a lace diagram $\bw = (w_1,\dots,w_n)$ for the dimension
vector $e$, we let $\Schub(w_1,\dots,w_n)$ be the product of the
Schubert polynomials $\Schub_{w_a}(x^a;x^{a-1})$ for all rightward
arrows $a$, as well as the polynomials $\Schub_{w_a}(\wt x^{a-1};\wt
x^a)$ for all leftward arrows $a$.
\begin{equation} \label{E:schublace}
  \Schub(w_1,\dots,w_n) =
  \left( \prod_{a: \delta(a)=1} \Schub_{w_a}(x^a;x^{a-1}) \right) \cdot
  \left( \prod_{a: \delta(a)=-1}\Schub_{w_a}(\wt x^{a-1}; \wt x^a) \right)
\end{equation}
Since each permutation $w_a$ is a partial permutation from $e_{h(a)}$
elements to $e_{t(a)}$ elements, it follows that the corresponding
Schubert polynomial receives the required number of variables.
Finally, for any $G$-orbit $\mu \subset V$ we define the polynomial
\[ Q_\mu = \sum_{(w_1,\dots,w_n)} \Schub(w_1,\dots,w_n) \]
where the sum is over all minimal lace diagrams for $\mu$.  Our main
result is the following theorem, which generalizes the equioriented
component formula proved in \cite{knutson.miller.ea:four} and
\cite{buch.feher.ea:positivity}.

\begin{thm} \label{T:comp}
  The polynomial $Q_\mu$ represents the $G$-equivariant cohomology
  class of the orbit closure $\ov{\mu}$ in $H^*_G(V)$.
\end{thm}

M.~Shimozono reports that he had speculated that this formula was
true, but had not been able to prove it.

\subsection{Degeneracy loci}

Theorem~\ref{T:comp} can be interpreted as a formula for degeneracy
loci defined by a quiver $F_\bull$ of vector bundle morphisms over a
non-singular complex variety $X$.  This quiver consists of a vector
bundle $F_i$ of rank $e_i$ for each vertex $i \in Q_0$, and a bundle
map $F_{t(a)} \to F_{h(a)}$ for each arrow $a \in Q_1$.  These bundle
maps define a section $s : X \to H$ to the bundle $\pi : H =
\bigoplus_{a \in Q_1} \Hom(F_{t(a)}, F_{h(a)}) \to X$.  Since each
fiber $\pi^{-1}(x)$ of $H$ is identical to the representation space
$V$, a $G$-orbit $\mu \subset V$ defines a Zariski closed subset
$H_\mu$ in $H$ as the union of the orbit closures $\ov\mu \subset V =
\pi^{-1}(x)$ for all $x \in X$.  The corresponding degeneracy locus in
$X$ is defined as the scheme theoretic inverse image $X_\mu =
s^{-1}(H_\mu)$.  We assume that the bundle maps of $F_\bull$ are
sufficiently generic, so that $X_\mu$ obtains its maximal possible
codimension $d(\mu)$ in $X$.

It follows from the definition of equivariant cohomology that the
cohomology class $[H_\mu] \in H^*(H)$ is given by the polynomial
$Q_\mu$, when the Chern roots of $\pi^* F_i$ are substituted for the
variables $x^i$.  Using Bobi{\'n}ski and Zwara's result that the orbit
closure $\ov\mu$ (and therefore $H_\mu$) is Cohen-Macaulay
\cite{bobinski.zwara:normality}, it follows from
\cite[Prop.~7.1]{fulton:intersection} that $[X_\mu] = s^*[H_\mu]$ in
$H^*(X)$, so the cohomology class of $X_\mu$ is also given by $Q_\mu$
when the variables $x^i$ are identified with the Chern roots of $F_i$.

If $X$ admits an ample line bundle $L$, then this formula remains true
in the Chow group of $X$.  In fact, by twisting the bundles $F_i$ with
a power of $L$, we may assume that these bundles are globally
generated.  In this case one can construct a bundle $Y = \bigoplus_{a
  \in Q_1} \Hom(B_{t(a)}, B_{h(a)})$ over a product of Grassmannians
$\prod_{i \in Q_0} \Gr^{e_i}(\C^N)$ with tautological quotient bundles
$B_i$ for $i \in Q_0$, such that the quiver $F_\bull$ on $X$ is the
pullback of the universal quiver $B_\bull$ on $Y$ along a morphism of
varieties $f : X \to Y$.  Since the Chow cohomology of $Y$ agrees with
singular cohomology, our formula for the Chow class of $X_\mu$ follows
from the identity $[X_\mu] = f^*[Y_\mu]$, which again uses that
$Y_\mu$ is Cohen-Macaulay.

\subsection{Symmetry of the component formula}

In order to apply the interpolation method from
\cite{feher.rimanyi:calculation} to prove Theorem~\ref{T:comp}, we
first need to show that the polynomial $Q_\mu$ belongs to the subring
$H^*_G(V)$ of symmetric polynomials in $H^*_T(V)$ (of course, this is
{\em implied\/} by Theorem~\ref{T:comp}). We prove this as in
\cite{buch.feher.ea:positivity}, except that there are more cases to
consider.

\begin{lemma} \label{L:sym}
  The polynomial $Q_\mu$ is separately symmetric in each set of
  variables $x^i$, $0\leq i\leq n$.
\end{lemma}
\begin{proof}
  We must show that for any $0 \leq i \leq n$ and $1 \leq j < e_i$,
  the divided difference operator $\partial^i_j =
  \partial_{x^i_j,x^i_{j+1}}$ maps $Q_\mu$ to zero.  We verify this
  using the identities (\ref{E:schubdef}) and (\ref{E:dualdef}) of
  Schubert polynomials.  Let $\bw = (w_1,\dots,w_n)$ be a minimal lace
  diagram for $\mu$.  For convenience, we identify each variable
  $x^i_k$ with dot $k$ from the top of column $i$.  Notice that if two
  line segments connected to $x^i_j$ and $x^i_{j+1}$ cross each other,
  then the minimality of the lace diagram implies that $0 < i < n$,
  and only the connections on one side of these dots are allowed to
  cross.
  
  Assume first that the line segments connecting $x^i_j$ and
  $x^i_{j+1}$ to dots of column $i-1$ cross each other.  Let $\bu =
  (u_1,\dots,u_n)$ be the lace diagram obtained from $\bw$ by removing
  this crossing.  In other words, we set $u_p = w_p$ for $p \neq i$,
  while $u_i = w_i s_j$ if arrow $i$ points right and $u_i = s_{e_i-j}
  w_i$ if arrow $i$ points left.  We claim that
\[ \partial^i_j(\Schub(w_1,\dots,w_n)) = \Schub(u_1,\dots,u_n) \,. \]
By using the identity $\partial^i_j(fg) = \partial^i_j(f)g$, which
holds for polynomials $f$ and $g$ such that $g$ is symmetric in
$\{x^i_j,x^i_{j+1}\}$, we need only check that $\partial^i_j$ maps the
$i$th factor of $\Schub(\bw)$ to the $i$th factor of $\Schub(\bu)$.
This follows from (\ref{E:schubdef}) when arrow $i$ points right and from
(\ref{E:dualdef}) when arrow $i$ points left.

One checks similarly that, if the line segments connecting $x^i_j$ and
$x^i_{j+1}$ to dots of column $i+1$ cross each other, then
$\partial^i_j(\Schub(\bw)) = - \Schub(\bu)$, where the lace diagram
$\bu$ is obtained from $\bw$ by removing this crossing.  Furthermore,
if none of the lines connected to $x^i_j$ and $x^i_{j+1}$ cross each
other, then $\partial^i_j(\Schub(\bw)) = 0$.

For each minimal lace diagram $\bw$ for $\mu$ in which the connections
to $x^i_j$ and $x^i_{j+1}$ from one side cross each other, one can
construct another minimal lace diagram $\bw'$ for $\mu$ by moving the
crossing to the opposite side of these dots.  The lemma follows from
this because $\partial^i_j(\Schub(\bw) + \Schub(\bw')) = 0$.
\end{proof}

\subsection{Existence of minimal lace diagrams}

The orbit-preserving transformation of lace diagrams exploited in the
proof of Lemma~\ref{L:sym} is illustrated by the following picture (of
parts of the extended lace diagrams):
\begin{equation} \label{E:genmove}
\raisebox{-9pt}{\pic{50}{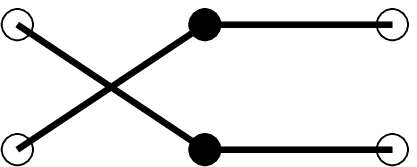}}
\ \ \  \longleftrightarrow \ \ \ 
\raisebox{-9pt}{\pic{50}{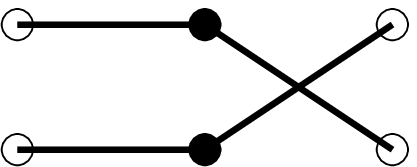}}
\end{equation}
These transformations played a similar role in
\cite{buch.feher.ea:positivity}.  Notice that the transformation
(\ref{E:genmove}) can be applied to any lace diagram, as long as the
middle dots and at least one from each column of outer dots are not in
the extended part of the diagram.

\begin{prop} \label{P:minlace}
  Let $\mu \subset V$ be any $G$-orbit.  Then there exists at least
  one minimal lace diagram representing $\mu$, and every minimal lace
  diagram for $\mu$ can be obtained from any other such diagram by
  using the transformations (\ref{E:genmove}).
\end{prop}
\begin{proof}
  Given any minimal lace diagram for $\mu$, we can use the
  transformations (\ref{E:genmove}) repeatedly, in left to right
  direction, until all crossings of the lace diagram involve the
  right hand side extension of one of the crossing strands.  It is
  therefore enough to prove that each orbit $\mu$ has a unique minimal
  lace diagram with this property.

  We will say that two (non-extended) strands {\em overlap\/} if both
  contain a dot in the same column.  Notice that if all crossings of a
  lace diagram occur in the extended part of the diagram, then the
  lace diagram is uniquely determined by specifying, for each pair of
  overlapping strands, which strand is placed above the other.  The
  uniqueness therefore follows from the observation that, if all
  crossings between two overlapping strands involve the right side
  extension of one of them, then this condition dictates which strand
  is over the other.

  Finally, to prove that a minimal lace diagram exists, it is
  sufficient to give a total order on the set of all pairs of integers
  $(i,j)$ with $0 \leq i\leq j \leq n$, such that if $(i,j) < (p,q)$
  and a strand from column $i$ to column $j$ is placed above a strand
  from column $p$ to column $q$, then these strands cross at most
  once, and if they do, the crossing must occur at the right side
  extension of one of them.  Such an ordering can be defined
  explicitly by writing $(i,j) < (p,q)$ if and only if one of the
  following conditions hold:
\begin{enumerate}
  \item $\delta(i) = -1$ and $\delta(p) = 1$.
  \item $\delta(i) = \delta(p) = -1$ and $i > p$.
  \item $\delta(i) = \delta(p) = 1$ and $i < p$.
  \item $i=p$ and $\delta(j+1)=-1$ and $\delta(q+1)=1$.
  \item $i=p$ and $\delta(j+1)=\delta(q+1)=-1$ and $j<q$.
  \item $i=p$ and $\delta(j+1)=\delta(q+1)=1$ and $j>q$.
\end{enumerate}
The following is an example of a minimal lace diagram where the
strands are arranged according to this order.
\[ \pic{37}{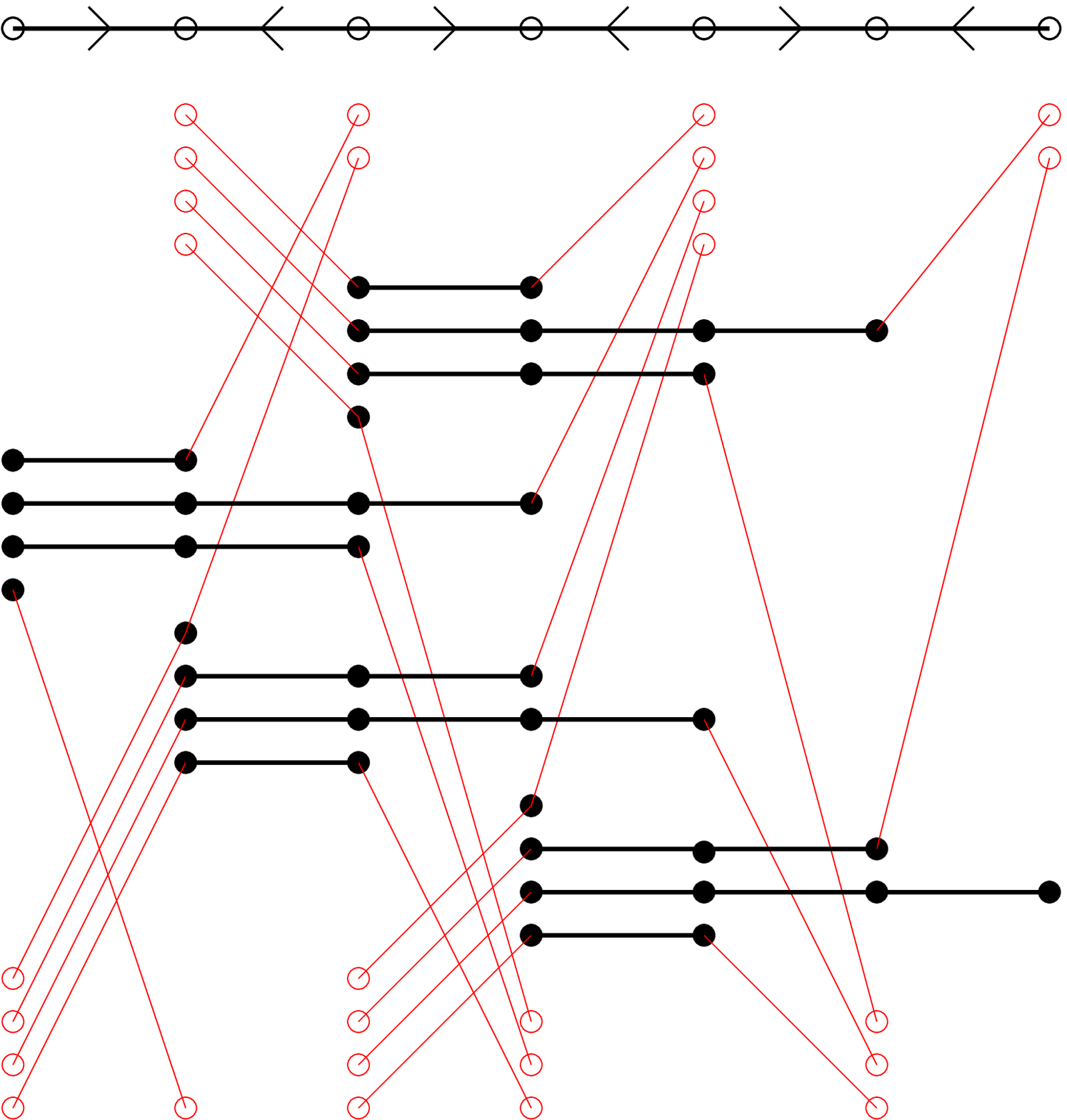} \vspace{-3mm} \]
\end{proof}

\section{Proof of the main theorem}
\label{S:proof}

\subsection{The interpolation method}

For each $G$-orbit $\mu \subset V$ we let $G_\mu$ denote the
stabilizer subgroup of a point $p_\mu$ in $\mu$.  The inclusion $G_\mu
\subset G$ induces a map $BG_\mu \to BG$, which gives an equivariant
restriction map $\phi_\mu : H^*_G(V) = H^*(BG) \to H^*(BG_\mu) =
H^*_G(\mu)$.  The {\em Euler class\/} $\euler(\mu) \in H_G^*(\mu)$ is
the top equivariant Chern class of the normal bundle to $\mu$ in $V$.
We will prove our formula for the class of $\ov \mu$ as an application
of the {\em interpolation method\/} of Feh\'er and Rim\'anyi.  This
method works more generally when $G$ is an arbitrary complex Lie group
acting on a vector space $V$ with finitely many orbits, such that
$\euler(\mu)$ is not a zero-divisor in $H^*_G(\mu)$ for each orbit
$\mu$.  We need the following statement
\cite[Thm.~3.5]{feher.rimanyi:calculation}.

\begin{thm} \label{T:interpolation}
  Let $\mu \subset V$ be a $G$-orbit.  The $G$-equivariant cohomology
  class of the closure of $\mu$ is the unique class $[\ov \mu] \in
  H^*_G(V)$ satisfying $\phi_\mu([\ov \mu]) = \euler(\mu) \in
  H^*_G(\mu)$ and $\phi_\eta([\ov \mu]) = 0$ for every $G$-orbit $\eta
  \subset V$ for which $\eta \neq \mu$ and $\codim \eta \leq \codim
  \mu$.
\end{thm}

\subsection{Description of the Euler class}

We also need a description of the restriction maps $\phi_\mu :
H^*_G(V) \to H^*_G(\mu)$ and Euler classes $\euler(\mu) \in
H^*_G(\mu)$ which was proved in \cite{feher.rimanyi:classes} for any
quiver of Dynkin type.  Fix a lace diagram $\bw$ representing the
orbit $\mu$, and choose variables $b_1,\dots,b_k$ corresponding to the
strands of $\bw$.  Then $H^*_G(\mu)$ can be identified with a subring
of the polynomial ring $\Z[b_1,\dots,b_k]$; this was done in
\cite[\S3]{feher.rimanyi:classes} by showing that $G_\mu$ has a
maximal torus of dimension equal to the number of strands.  By
\cite[Prop.~3.10]{feher.rimanyi:classes}, the restriction map
$\phi_\mu : H^*_G(V) \to H^*_G(\mu)$ extends to a ring homomorphism
$\phi_{\bw} : H^*_T(V) \to \Z[b_1,\dots,b_k]$, which maps $x^i_j$ to
the variable of the strand passing through dot $j$ from the top of
column $i$ of the lace diagram.  This map $\phi_\bw$ depends on the
chosen lace diagram $\bw$ for $\mu$.

To describe the Euler class $\euler(\mu) \in H^*_G(\mu)$, we need some
definitions for quiver representations with arbitrary dimension
vectors.  Let $\phi = (\phi_1,\dots,\phi_n)$ and $\phi' =
(\phi'_1,\dots,\phi'_n)$ be representations of $Q$ with dimension
vectors $e = (e_0,\dots,e_n)$ and $e' = (e'_0,\dots,e'_n)$.  A
homomorphism $\alpha : \phi \to \phi'$ is a tuple $\alpha =
(\alpha_0,\dots,\alpha_n)$ of linear maps $\alpha_i : \C^{e_i} \to
\C^{e'_i}$ such that $\alpha_{h(a)} \phi_a = \phi'_a \alpha_{t(a)}$
for all arrows $a$.  The set $\Hom(\phi,\phi')$ of all such
homomorphisms is a complex vector space.  By using an injective
resolution of the representation $\phi'$, one can also define the
extension module $\Ext(\phi,\phi') = \Ext^1(\phi,\phi')$.  Let $E_Q$
be the Euler form defined by $E_Q(\phi,\phi') = E_Q(e,e') = \sum_{i
  \in Q_0} e_i e'_i - \sum_{a \in Q_1} e_{t(a)} e'_{h(a)}$.  The
homomorphism and extension modules are related by the identity
\cite{ringel:representations}
\begin{equation} \label{E:eform}
  E_Q(\phi,\phi') = \dim \Hom(\phi,\phi') - \dim \Ext(\phi,\phi') \,.
\end{equation}

A quiver representation is {\em indecomposable\/} if it cannot be
written as a direct sum of other quiver representations.  For a quiver
of Dynkin type, the indecomposable representations correspond to the
positive roots of the corresponding root system
\cite{gabriel:unzerlegbare} (see also
\cite{bernstein.gelfand.ea:coxeter}).  For our quiver of type $A$,
there is one indecomposable representation $X^{ij}$ for each pair of
integers $(i,j)$ with $0 \leq i \leq j \leq n$.  The dimension vector
of $X^{ij}$ assigns the dimension $1$ to all vertices $k \in Q_0$ with
$i \leq k \leq j$, and assigns dimension zero to all other vertices.
For each arrow $i < a \leq j$, the map $X^{ij}_a : \C \to \C$ is the
identity.  Given a $G$-orbit $\mu \subset V$, the indecomposable
summands in the decomposition of a representation $\phi\in\mu$
correspond to the strands in the lace diagram for $\mu$.  More
canonically, the multiplicity of $X^{ij}$ in $\phi$ is equal to the
number of strands $s_{ij}(\mu)$ from column $i$ to column $j$.

We can now state the formula for the Euler class $\euler(\mu)$, using
the above described embedding $H^*_G(\mu) \subset \Z[b_1,\dots,b_k]$.
For each pair of variables $b_p, b_q$ we let $\Ext(b_p,b_q)$ denote
the extension module of the indecomposable representations
corresponding to the strands of $b_p$ and $b_q$.  The following was
proved in \cite[Cor.~3.13]{feher.rimanyi:classes}.

\begin{prop} \label{P:euler}
The Euler class of the $G$-orbit $\mu \subset V$ is given by
\[ \euler(\mu) = \prod_{1 \leq p,q \leq k}
   (b_p - b_q)^{\dim \Ext(b_q,b_p)} \,.
\]
\end{prop}

\subsection{Proof of the formula}

We need to compute the dimension of an extension module
$\Ext(X^{ij},X^{pq})$.  Let $N(X^{ij},X^{pq})$ denote the number of
arrows $a \in Q_1$ such that $t(a) \in [i,j]$, $h(a) \in [p,q]$, and
such that $h(a) \not \in [i,j]$ or $t(a) \not \in [p,q]$.

\begin{lemma} The dimension of the extension module of the
  indecomposable representations $X^{ij}$ and $X^{pq}$ is given by
\[ \dim \Ext(X^{ij},X^{pq}) = \begin{cases}
  1 & \text{if $[i,j]\cap [p,q] \neq \emptyset$ and $N(X^{ij},X^{pq})
  = 2$} \\
  1 & \text{if $[i,j]\cap [p,q] = \emptyset$ and $N(X^{ij},X^{pq}) =
  1$} \\
  0 & \text{otherwise.}
  \end{cases}
\]
\end{lemma}
\begin{proof}
  If $[i,j] \cap [p,q] = \emptyset$ then $\Hom(X^{ij},X^{pq}) = 0$ and
  $N(X^{ij},X^{pq}) \leq 1$.  The lemma follows from (\ref{E:eform})
  because $E_Q(X^{ij},X^{pq}) = - N(X^{ij},X^{pq})$.

  Otherwise $[i,j] \cap [p,q] \neq \emptyset$, in which case we have
  $N(X^{ij},X^{pq}) \leq 2$.  It follows from the definition that
  $\Hom(X^{ij},X^{pq}) = \C$ if $N(X^{ij},X^{pq}) = 0$, while
  $\Hom(X^{ij},X^{pq}) = 0$ otherwise.  The lemma now follows because
  $E_Q(X^{ij},X^{pq}) = 1 - N(X^{ij},X^{pq})$.
\end{proof}

Another way to state this lemma is that $\Ext(X^{ij},X^{pq})$ is
non-zero (with dimension one) exactly when two strands corresponding
to $X^{ij}$ and $X^{pq}$ are forced to cross each other (see the proof
of Lemma~\ref{L:mincross}), and when $(i,j) < (p,q)$ in the order used
in the proof of Proposition~\ref{P:minlace}.  Notice also that if two
such strands have a single crossing point, then the slope at the
crossing point of the strand corresponding to $X^{pq}$ is larger than
the slope of the strand corresponding to $X^{ij}$.  As a consequence
we obtain the following description of the Euler class $\euler(\mu)$.

\begin{cor} \label{C:euler}
  Let $\bw$ be a minimal lace diagram for the $G$-orbit $\mu \subset
  V$, and let $H^*_G(\mu) \subset \Z[b_1,\dots,b_k]$ be the
  corresponding inclusion of rings.  Then the Euler class $\euler(\mu)
  \in H^*_G(\mu)$ is the product of all factors $(b_p - b_q)$ for
  which the strands of $b_p$ and $b_q$ cross each other and the strand
  of $b_p$ has the highest slope at the crossing point.
\end{cor}

\begin{cor} \label{C:codim}
  The codimension of the $G$-orbit $\mu \subset V$ is equal
  to the length $d(\mu)$ of any minimal lace diagram for $\mu$.
\end{cor}

\begin{proof}[Proof of Theorem~\ref{T:comp}]
  It follows from Lemma~\ref{L:sym} that $Q_\mu$ is an element of
  $H^*_G(V)$.  According to Theorem~\ref{T:interpolation}, we need to
  prove that $\phi_\mu(Q_\mu) = \euler(\mu)$ and that
  $\phi_\eta(Q_\mu) = 0$ for any $G$-orbit $\eta \subset V$ such that
  $\eta \neq \mu$ and $\codim \eta \leq \codim \mu$.  It is enough to
  show that if $\bu$ is a minimal lace diagram for $\mu$ and $\bw$ is
  any lace diagram for the same dimension vector $e$ such that
  $\ell(\bw) \leq \ell(\bu)$, then we have
\[ \phi_\bw(\Schub(\bu)) = \begin{cases}
  \euler(\mu) & \text{if $\bw = \bu$} \\
  0 & \text{if $\bw \neq \bu$.}
\end{cases}\]
For any arrow $a \in Q_1$, $\phi_\bw$ maps the $a$th factor of
$\Schub(\bu)$ to the specialized Schubert polynomial
$\Schub_{u_a}(b_{w_a(1)},\dots,b_{w_a(m)}; b_1,\dots,b_m)$, where
$b_1,\dots,b_m$ denote the variables of strands of $\bw$ connecting
column $a-1$ to column $a$.  If $\delta(a) = 1$ then $b_1,\dots,b_m$
correspond to the strands passing through the dots of column $a-1$
ordered from top to bottom, and starting with the first non-extended
dot.  If $\delta(a) = -1$, then we use the strands passing through the
dots of column $a$ in bottom to top order, starting with the lowest
non-extended dot.  Now the specialization $\Schub_{u_a}(b_{w_a};b)$ is
zero unless $u_a \leq w_a$ in the Bruhat order.  Since $\ell(\bw) \leq
\ell(\bu)$, it follows that $\phi_\bw(\Schub(\bu))$ is zero unless
$\bu = \bw$.  Furthermore, it follows from (\ref{E:spec}) that
$\Schub_{u_a}(b_{u_a}; b)$ is equal to the product of the factors
$(b_p - b_q)$ of Corollary~\ref{C:euler} for which the strands of
$b_p$ and $b_q$ cross between column $a-1$ and column $a$.  This shows
that $\phi_\bu(\Schub(\bu)) = \euler(\mu)$, and finishes the proof.
\end{proof}

\begin{remark}
  In most of our pictures of lace diagrams, the columns of
  (non-extended) dots have been aligned at the top.  However, the
  definition of extended lace diagrams makes it natural to
  bottom-align two consecutive columns if they are connected by a
  leftward arrow, while columns connected with a rightward arrow are
  top-aligned as usual.  When this convention is used, a minimal lace
  diagram for the open orbit in the representation space $V$ can be
  obtained by simply drawing all possible horizontal lines between
  dots of consecutive columns.  For example, the open orbit for the
  quiver $(\circ \to \circ \ot \circ \to \circ \to \circ \ot \circ \to
  \circ \ot \circ)$ with dimension vector $e = (2,4,3,2,3,4,2,3)$ is
  represented by the following minimal lace diagram.
\[ \pic{40}{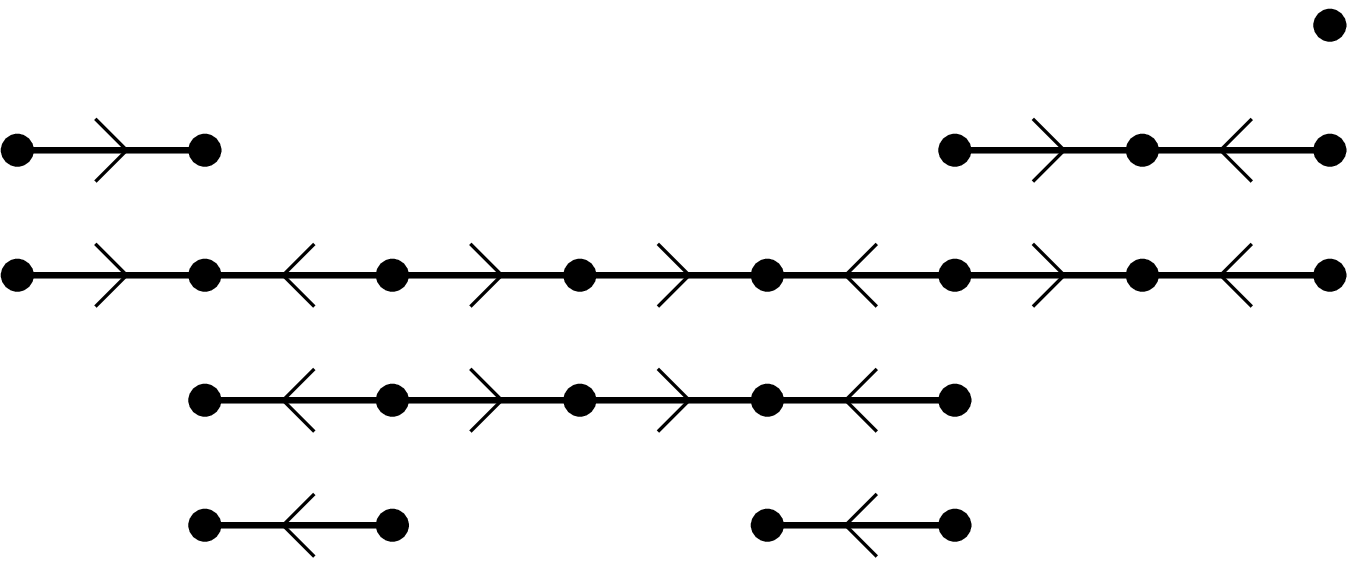} \]
\end{remark}

\section{A conjectural $K$-theoretic component formula}
\label{S:ktheory}

For an equioriented quiver of type $A$, Buch has proved a
$K$-theoretic generalization of the component formula, which expresses
the equivariant Grothen\-dieck class of an orbit closure as an
alternating sum of products of Grothendieck polynomials
\cite{buch:alternating}.  This sum is over all {\em
  KMS-factorizations\/} for the orbit, which generalize the
equioriented minimal lace diagrams from \cite{knutson.miller.ea:four}.
A limit of the $K$-theoretic component formula in terms of stable
Grothendieck polynomials was also obtained by Miller
\cite{miller:alternating}.  It was proved in
\cite{buch.feher.ea:positivity} that all KMS-factorizations for an
orbit can be obtained from a minimal lace diagram for the orbit by
applying a series of transformations:
\[
\raisebox{-9pt}{\pic{50}{tx_.eps}} \ \ \ \longleftrightarrow \ \ \ 
\raisebox{-9pt}{\pic{50}{t_x.eps}} \ \ \ \longleftrightarrow \ \ \ 
\raisebox{-9pt}{\pic{50}{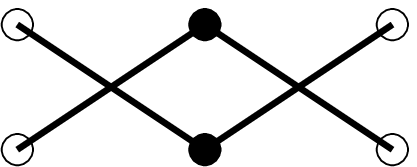}}
\]
In these transformations, the middle dots and at least one of the
outer dots from each side must be outside the extended part of the
diagram, and furthermore the two middle dots must be consecutive in
their column.

Given an orbit $\mu \subset V$ of representations of an arbitrary
quiver of type $A$, define a {\em $K$-theoretic lace diagram\/} for
this orbit to be any lace diagram that can be obtained from a minimal
lace diagram representing $\mu$ by using these transformations.  For
such a diagram $\bw$, we define a Laurent polynomial $\Groth(\bw)$ by
the expression (\ref{E:schublace}), except that each Schubert
polynomial $\Schub_w(x;y)$ is replaced with the Grothendieck (Laurent)
polynomial $\Groth_w(x;y)$ from
\cite{lascoux.schutzenberger:structure}.  This polynomial is defined
by the recursive identities $(x_i-x_{i+1}) \Groth_w(x;y) = x_i
\Groth_{w s_i}(x;y) - x_{i+1} \Groth_{w s_i}(x_{s_i}; y)$ when $w(i) <
w(i+1)$, as well as the expression $\Groth_{w_0}(x;y) = \prod_{i+j
  \leq m} (1-y_i/x_j)$ when $w_0$ is the longest permutation in $S_m$.
Given a maximal torus $T \subset G$, we identify the variable $x^i_j$
of $\Groth(\bw)$ with the $T$-equivariant class of a line bundle $V
\times \C \to V$ with the action of $T$ given by $t . (\phi,z) =
(t.\phi,t^i_j z)$, where $(t^i_1,\dots,t^i_{e_i})$ are chosen
coordinates on $T \cap \GL(E_i)$.  We finish this paper by posing the
following conjecture, which generalizes Theorem~\ref{T:comp} as well
as \cite[Thm.~6.3]{buch:alternating}.

\begin{conj}
  The $T$-equivariant Grothendieck class of $\ov \mu$ is given by
\[ [{\mathcal O}_{\ov \mu}] =
   \sum_{\bw} (-1)^{\ell(\bw)-d(\mu)} \Groth(\bw)
\]
where the sum is over all $K$-theoretic lace diagrams for $\mu$.
\end{conj}

Aside from the analogies with known identities, this conjecture is
motivated by the fact that one is naturally led to the $K$-theoretic
transformations when attempting to prove that a linear combination of
the products $\Groth(\bw)$ of Grothendieck polynomials is symmetric in
each set of variables $x^i$.  We note that it is possible to compute
the $T$-equivariant Grothendieck class of $\ov\mu$ by using a
resolution of singularities of this locus like in
\cite{buch.fulton:chern, buch:grothendieck}, but there is no known way
to derive the conjectured formula from this approach, even for
equioriented quivers.


\providecommand{\bysame}{\leavevmode\hbox to3em{\hrulefill}\thinspace}
\providecommand{\MR}{\relax\ifhmode\unskip\space\fi MR }
\providecommand{\MRhref}[2]{%
  \href{http://www.ams.org/mathscinet-getitem?mr=#1}{#2}
}
\providecommand{\href}[2]{#2}


\end{document}